

\hsize=16.4cm
\vsize=22.3cm
\hoffset=-0.1cm
\voffset=0.5cm

\nopagenumbers

\input amssym.def
\input amssym.tex

\font\teneufm=eufm10
\font\seveneufm=eufm7
\font\fiveeufm=eufm5 \newfam\eufmfam
\def\eufm{\fam\eufmfam\teneufm}
\textfont\eufmfam=\teneufm
\scriptfont\eufmfam=\seveneufm 
\scriptscriptfont\eufmfam=\fiveeufm

\def\Re{{\rm Re}\,}

\def\scr#1{{\scriptstyle{#1}}}

\def\B#1{{\Bbb #1}}
\def\ca#1{{\cal #1}}

\font\title=cmbx12 at 14pt
\font\author=cmr12
\font\srm=cmr8 at 9pt
 at11pt
\font\small=cmr7

\def\rightheadline{\hfil{\srm}\hfil\folio}
\def\leftheadline{\rm\folio\hfil{\srm}\hfil}
\def\emptyheadline{\hfil{}}
\headline{\rm\ifnum\pageno=1 \emptyheadline\else
\ifodd\pageno \rightheadline \else \leftheadline\fi\fi}

\def\firstpage{\hss{\vbox to 1cm{\vfil\hbox{\rm\folio}}}\hss}
\def\emptyfootline{\hfil}
\footline{\ifnum\pageno=1\firstpage\else
\emptyfootline\fi}

\centerline{\title Small Gaps between Primes Exist}
\vskip 1cm
\centerline{By D.A. Goldston, Y. Motohashi,
J. Pintz, and C.Y. Y{\i}ld{\i}r{\i}m}
\footnote{}{\small
\hskip -0.7cm
1991 Mathematics Subject Classification: Primary 11N05; Secondary 11P32
\par
\hskip -0.7cm
Key Words and phrases: prime number
\par
\hskip -0.7cm
The first author was supported by NSF grant DMS-0300563, the NSF Focused Research 
Group grant 0244660, and the American Institute of Mathematics;
the second author by KAKENHI 15540047 and Nihon University Research grant; 
the third author by OTKA grants No.\ T38396, T43623, T49693 and
Balaton program; the fourth author by T\"UBITAK}
\vskip 1cm
{\bf Abstract.} In the recent preprint [3], Goldston, Pintz, and Y{\i}ld{\i}r{\i}m 
established, among other things, $$ \liminf_{n\to\infty}{p_{n+1}-p_n\over\log 
p_n}=0,\leqno(0) $$ with $p_n$ the $n$th prime. In the present article, which is 
essentially self-contained, we shall develop a simplified account of the method used 
in [3]. While [3] also includes quantitative versions of $(0)$, we are concerned here 
solely with proving the qualitative $(0)$, which still exhibits all the essentials 
of the method. We also show here that an improvement of the Bombieri--Vinogradov prime 
number theorem would give rise infinitely often to bounded differences between 
consecutive primes. We include  a short expository last section.  Detailed 
discussions of quantitative  results and a historical review will appear in the 
publication version of [3] and its continuations.
\vskip 0.6cm
\centerline{\bf 1. Basic Lemma}
\medskip
In this section we shall prove an asymptotic formula relevant to 
Selberg's sieve, which is to be modified so as to involve primes
in the next section. The two asymptotic formulas thus obtained will be
combined in a simple weighted sieve setting, and give rise to $(0)$ in
the third section. 
\smallskip
Let $N$ be a parameter increasing monotonically to infinity. There
are four other basic parameters $H, R, k, \ell$ in our discussion. We
impose the following conditions to them:
$$
H\ll\log N\ll\log R\le\log N, \leqno(1.1)
$$
and
$$
\hbox{integers $k, \ell> 0$ are arbitrary but bounded.}
\leqno(1.2)
$$
To prove a quantitative assertion superseding $(0)$, we need to regard
$k,\ell$ as functions of $N$; but for our present
purpose  the above is sufficient. All implicit constants in the sequel
are possibly dependent on $k,\ell$ at most; and besides, the symbol $c$
stands for a positive constant with the same dependency, whose value may
differ at each occurrence.
\medskip
Let
$$
\ca{H}=\{h_1,h_2,\ldots,h_k\}\subseteq[1,H]\cap\B{Z},\leqno(1.3)
$$
with $h_i\ne h_j$ for $i\ne j$. Let us put, for a prime $p$,
$$
\Omega(p)=\{\hbox{different residue classes among $-h(\bmod\,p)$,
$h\in\ca{H}$}\}
\leqno(1.4)
$$
and write $n\in\Omega(p)$ instead of $n\,(\bmod\,p)\in\Omega(p)$. We call  
$\ca{H}$  admissible if
$$
|\Omega(p)|<p\quad \hbox{for all $p$},\leqno(1.5)
$$
and assume this unless otherwise stated.
We extend $\Omega$ multiplicatively, so that 
$n\in\Omega(d)$ with square-free $d$ if and only if
$n\in\Omega(p)$ for all $p|d$, which is equivalent to
$$
d| P(n;\ca{H}),\quad P(n;\ca{H})=(n+h_1)(n+h_2)\cdots(n+h_k).
\leqno(1.6)
$$
Also, we put, with $\mu$ the M\"obius function,
$$
\lambda_R(d;a)=\cases{\hfil 0 & if $d>R$,\cr
\displaystyle{1\over a!}\mu(d)
(\log R/d)^a & if $d\le R$,}\leqno(1.7)
$$
and
$$
\Lambda_R(n;\ca{H},a)
=\sum_{n\in\Omega(d)}\lambda_R(d;a)= {1\over a!}\sum_{\scriptstyle  d| 
P(n;\ca{H}) \atop \scriptstyle d\leq R}\mu(d)
(\log R/d)^a .\leqno(1.8)
$$
\par
With this, we shall consider the evaluation of
$$
\sum_{N<n\le 2N}\Lambda_R(n;\ca{H},k+\ell)^2.\leqno(1.9)
$$
Expanding out the square, we have
$$
\sum_{d_1,d_2}\lambda_R(d_1;k+\ell)\lambda_R(d_2;k+\ell)
\sum_{\scr{N<n\le 2N}\atop\scr{n\in\Omega(d_1),n\in\Omega(d_2)}}1.
\leqno(1.10)
$$
The condition on $n$ is equivalent to
$n\in\Omega([d_1,d_2])$, with $[d_1,d_2]$ the least common multiple of
the two integers; and 
$$
\sum_{N<n\le 2N}\Lambda_R(n;\ca{H},k+\ell)^2=
N\ca{T}+O\left(\Big(\sum_{d}|\Omega(d)||\lambda_R(d;k+\ell)|
\Big)^2\right),
\leqno(1.11)
$$
in which
$$
\ca{T}=\sum_{d_1,d_2}{|\Omega([d_1,d_2])|\over[d_1,d_2]}
\lambda_R(d_1;k+\ell)\lambda_R(d_2;k+\ell).\leqno(1.12)
$$
We have $|\Omega(d)|\le\tau_k(d)$
with the divisor function $\tau_k$. Thus
$$
\sum_{N<n\le 2N}\Lambda_R(n;\ca{H},k+\ell)^2=
N\ca{T}+O\left(R^2(\log R)^c\right).
\leqno(1.13)
$$
\par
On noting that for $a\ge1$
$$
\lambda_R(d;a)={\mu(d)\over2\pi i}\int_{(1)}\left({R\over
d}\right)^s{ds\over s^{a+1}},\leqno(1.14)
$$
with $(\alpha)$ the vertical line in the complex plane passing through 
$\alpha$, we have
$$
\ca{T}={1\over(2\pi i)^2}\int_{(1)}\!\int_{(1)}F(s_1,s_2;\Omega)
{R^{s_1+s_2}\over(s_1s_2)^{k+\ell+1}}ds_1ds_2,\leqno(1.15)
$$
where
$$
\leqalignno{
F(s_1,s_2;\Omega)&=\sum_{d_1,d_2}\mu(d_1)\mu(d_2)
{|\Omega([d_1,d_2])|\over[d_1,d_2]d_1^{s_1}d_2^{s_2}}&(1.16)\cr
&=\prod_p\left(1-{|\Omega(p)|\over p}\left({1\over p^{s_1}}+
{1\over p^{s_2}}-{1\over p^{s_1+s_2}}
\right)\right)
}
$$
 in the region of absolute convergence.
\par
Since $|\Omega(p)|=k$ for $p>H$, we put
$$
G(s_1,s_2;\Omega)=F(s_1,s_2;\Omega)\left(\zeta(s_1+1)\zeta(s_2+1)
\over\zeta(s_1+s_2+1)\right)^k,\leqno(1.17)
$$
with $\zeta$ the Riemann zeta-function.
This is regular and bounded for $\Re s_1,\Re s_2>-c$. 
In particular, we have the singular series
$$
{\eufm S}(\ca{H})=G(0,0;\Omega)=\prod_p\left(1-{|\Omega(p)|\over
p}\right)\left(1-{1\over p}\right)^{-k},\leqno(1.18)
$$
which does not vanish because of $(1.5)$. We have the
bound
$$
G(s_1,s_2;\Omega)\ll
\exp(c(\log N)^{-2\sigma}\log\log\log N),\leqno(1.19)
$$ 
with $\min(\Re s_1,\Re s_2,0)=\sigma\ge -c$,
as can be seen via the Euler product expansion of the 
right side of $(1.17)$. In fact, the part corresponding to those $p>H$ is
uniformly bounded in the indicated region since
$|\Omega(p)|=k$. For $k^2<p\le H$, the logarithm of each $p$-factor 
is estimated to be $\ll H^{-2\sigma}\sum_{p\le H}p^{-1}$;
and the treatment of those $p\le k^2$ is trivial.
Note that the restrictions $(1.1)$ and $(1.2)$ are relevant here.  
\medskip
Now, we have in place of $(1.15)$
$$
\ca{T}={1\over(2\pi i)^2}\int_{(1)}\!\int_{(1)}G(s_1,s_2;\Omega)
\left(\zeta(s_1+s_2+1)\over\zeta(s_1+1)\zeta(s_2+1)\right)^k
{R^{s_1+s_2}\over(s_1s_2)^{k+\ell+1}}ds_1ds_2.\leqno(1.20)
$$
Let us put $U=\exp(\sqrt{\log N})$, and shift the $s_1$ and
$s_2$-contours to the vertical lines $(\log U)^{-1}+it$ and to
$(2\log U))^{-1}+it$, $t\in\B{R}$, respectively.
We truncate them to $|t|\le U$ and
$|t|\le U/2$, and denote the results
by $L_1$ and $L_2$, respectively. On noting $(1.1)$ and $(1.19)$, we
have  readily that
$$
\leqalignno{
\ca{T}={1\over(2\pi i)^2}\int_{L_2}\!\int_{L_1}&G(s_1,s_2;\Omega)
\left(\zeta(s_1+s_2+1)\over\zeta(s_1+1)\zeta(s_2+1)\right)^k
{R^{s_1+s_2}\over(s_1s_2)^{k+\ell+1}}ds_1ds_2&(1.21)\cr
&+O\left(\exp(-c\sqrt{\log N})\right).
}
$$
We then shift the $s_1$-contour to $L_3: -(\log U)^{-1}+it$,
$|t|\le U$; necessary facts about the functions $\zeta$ and $1/\zeta$
can be found  in [4, p.\ 53] (or p.\ 60 in the Second Edition). We encounter
singularities at
$s_1=0$ and
$s_1=-s_2$, which are poles of orders
$\ell+1$ and $k$, respectively. We have
$$
\ca{T}={1\over2\pi i}\int_{L_2}\left\{\mathop{\rm Res}_{s_1=0}
+\mathop{\rm Res}_{s_1=-s_2}\right\}ds_2+
O\left(\exp(-c\sqrt{\log N})\right),\leqno(1.22)
$$
in which we have used $(1.19)$.
\par
We rewrite the residue, and have
$$
\mathop{\rm Res}_{s_1=-s_2}={1\over2\pi i}\int_{C(s_2)}
G(s_1,s_2;\Omega)
\left(\zeta(s_1+s_2+1)\over\zeta(s_1+1)\zeta(s_2+1)\right)^k
{R^{s_1+s_2}\over(s_1s_2)^{k+\ell+1}}ds_1,\leqno(1.23)
$$
with the circle $C(s_2): |s_1+s_2|=(\log N)^{-1}$. By $(1.19)$, we have
$G(s_1,s_2;\Omega)\ll(\log\log N)^c$; and $R^{s_1+s_2}\ll1$,
$\zeta(s_1+s_2+1)\ll\log N$. Also, since $|s_2|\ll|s_1|\ll |s_2|$, we
have
$(s_1\zeta(s_1+1))^{-1}\ll (|s_2|+1)^{-1}\log(|s_2|+2)$, 
{\it loc.cit\/}. Thus
$$
\mathop{\rm Res}_{s_1=-s_2}\ll(\log N)^{k-1}(\log\log N)^c
\left({\log (|s_2|+2)\over |s_2|+1}\right)^{2k}
|s_2|^{-2\ell-2}.\leqno(1.24)
$$
Inserting this into $(1.22)$, we get
$$
\ca{T}={1\over2\pi i}\int_{L_2}\left\{\mathop{\rm Res}_{s_1=0}
\right\}ds_2+
O\left((\log N)^{k+\ell}\right).\leqno(1.25)
$$
\par
To evaluate the last integral, we put 
$$
Z(s_1,s_2)=G(s_1,s_2;\Omega)\left({(s_1+s_2)\zeta(s_1+s_2+1)\over
s_1\zeta(s_1+1) s_2\zeta(s_2+1)}\right)^k,\leqno(1.26)
$$
which is regular in a neighborhood of the point $(0,0)$.
Then we have
$$
\mathop{\rm Res}_{s_1=0}
={R^{s_2}\over\ell!s_2^{\ell+1}}
\left({\partial\over\partial s_1}\right)^{\ell}_{s_1=0}
\left\{{Z(s_1,s_2)\over(s_1+s_2)^k}R^{s_1}\right\}.\leqno(1.27)
$$
We insert this into $(1.25)$ and shift the contour to $L_4:
-(2\log U)^{-1}+it$, $|t|\le U/2$. We see the new integral
is  $O(\exp(-c\sqrt{\log N}))$; the necessary bound for the 
integrand could be obtained in much the same way as $(1.24)$. Thus
$$
\leqalignno{
\ca{T}&=\mathop{\rm Res}_{s_2=0}\mathop{\rm Res}_{s_1=0}
+O((\log N)^{k+\ell})&(1.28)\cr
&={1\over(2\pi i)^2}
\int_{C_2}\int_{C_1}{Z(s_1,s_2)R^{s_1+s_2}\over
(s_1+s_2)^k(s_1s_2)^{\ell+1}}ds_1ds_2+O((\log N)^{k+\ell}),
}
$$
where $C_1$, $C_2$ are the circles $|s_1|=\rho$, 
$|s_2|=2\rho$, with a small $\rho>0$. We write
$s_1=s$, $s_2=s\xi$. Then  the double integral is equal to
$$
{1\over(2\pi i)^2}
\int_{C_3}\int_{C_1}{Z(s,s\xi)R^{s(\xi+1)}\over
(\xi+1)^k\xi^{\ell+1}s^{k+2\ell+1}}dsd\xi,\leqno(1.29)
$$
where $C_3$ is the circle $|\xi|=2$. This is
equal to
$$
{Z(0,0)\over2\pi i(k+2\ell)!}(\log R)^{k+2\ell}
\int_{C_3}{(\xi+1)^{2\ell}\over
\xi^{\ell+1}}d\xi+O((\log N)^{k+2\ell-1}(\log\log N)^c),\leqno(1.30)
$$
where we have used $(1.19)$.
\medskip
Hence, we have obtained our basic implement:
\medskip
\noindent
{\bf Lemma 1.} {\it Provided $(1.1)$, $(1.2)$, and $R\le N^{1/2}/(\log
N)^C$ hold with a sufficiently large $C>0$ depending only on $k$ and $\ell$, we have
$$
\sum_{N<n\le 2N}\Lambda_R(n;\ca{H},k+\ell)^2
={{\eufm S}(\ca{H})\over(k+2\ell)!}
{2\ell\choose\ell}N(\log R)^{k+2\ell}+O(N(\log N)^{k+2\ell-1}(\log\log
N)^c).\leqno(1.31)
$$
}
\vskip 1cm
\centerline{\bf 2. Twist with Primes} 
\medskip
Next, let us put
$$
\varpi(n)=\cases{ \log n & if $n$ is a prime,\cr
\hfil 0& otherwise;}\leqno(2.1)
$$
and consider the evaluation of the sum
$$
\sum_{N<n\le 2N}\varpi(n+h)\Lambda_R(n;\ca{H},k+\ell)^2,\leqno(2.2)
$$
with an arbitrary positive integer $h\le H$.  We observe that by $(1.6)$
this is  equal to 
$$
\sum_{N<n\le 2N}\varpi(n+h)
\Lambda_R(n;\ca{H}\!\setminus\!\{h\},k+\ell)^2,\leqno(2.3)
$$
provided $R<N$; in fact, if $\varpi(n+h)\ne0$ and $h\in\ca{H}$, then the
factor $n+h$ of $P(n;\ca{H})$ is irrelevant in computing
$\Lambda_R(n;\ca{H};k+\ell)$.
\medskip
We shall work on the
assumption: There exists an absolute constant
$0<\theta<1$ such that we have, for any fixed $A>0$,
$$
\qquad\sum_{q\le x^\theta}\max_{y\le x}\max_{\scr{a}\atop\scr{(a,q)=1}}
\left|\vartheta^*(y;a,q)-y/\varphi(q)\right|\ll x/(\log x)^A,
\quad\vartheta^*(y;a,q)=\sum_{\scr{y<n\le 2y}\atop
\scr{n\equiv a\bmod q}}\varpi(n),\leqno(2.4)
$$
with the implicit constant depending only on $A$. We assume that
$$
R\le N^{\theta/2}/(\log N)^A.\leqno(2.5)
$$
In particular, $(2.3)$ implies that we may assume also that
$h\notin\ca{H}$.
\medskip
With this, expanding out the square in $(2.2)$, we see that the sum is
equal to
$$
\sum_{d_1,d_2}\lambda_R(d_1;k+\ell)\lambda_R(d_2;k+\ell)
\sum_{b\in\Omega([d_1,d_2])}\delta((b+h,[d_1,d_2]))
\vartheta^*(N;b+h,[d_1,d_2]),\leqno(2.6)
$$
where $\delta(x)$ is the unit measure placed at $x=1$, because
$\vartheta^*(N,b+h,[d_1,d_2])=0$ if $b+h$ and $[d_1,d_2]$ are
not coprime.
Then, by $(2.4)$, this is equal to
$$
N\ca{T}^*+O(N/(\log N)^{A/3}),\leqno(2.7)
$$
with
$$
\ca{T}^*=\sum_{d_1,d_2}{\lambda_R(d_1;k+\ell)\lambda_R(d_2;k+\ell)\over
\varphi([d_1,d_2])}
\sum_{b\in\Omega([d_1,d_2])}\delta((b+h,[d_1,d_2])).\leqno(2.8)
$$
The error term in $(2.7)$ might require an explanation: We divide the
sum in $(2.6)$ into two parts according as 
$|\Omega([d_1,d_2])|=\tau_k([d_1,d_2])\le (\log N)^{A/2}$ and otherwise.
To the first part we apply $(2.4)$, while the second part is 
$$
\ll N(\log R)^{2(k+\ell)}\log N\sum_{d_1,d_2\le R}
{\tau_k([d_1,d_2])\over(\log N)^{A/2}}
{|\Omega([d_1,d_2])|\over[d_1,d_2]}\ll N/(\log N)^{A/3},\leqno(2.9)
$$
provided $A$ is sufficiently large.
\medskip
It remains for us to evaluate $\ca{T}^*$. The inner sum of $(2.8)$ is
equal to
$$
\prod_{p|[d_1,d_2]}\left(\sum_{b\in\Omega(p)}\delta((b+h,p))\right)
=\prod_{p|[d_1,d_2]}(|\Omega^+(p)|-1).\leqno(2.10)
$$
Here $\Omega^+$ corresponds to the set $\ca{H}^+=\ca{H}\cup\{h\}$.
In fact, $\delta((b+h,p))$ vanishes if and only if $-h\in\Omega(p)$;
and the latter is equivalent to $\Omega(p)=\Omega^+(p)$.
Note that the analogue of $(1.5)$ for
$\Omega^+$ could be violated. Nevertheless,  we have, as before,
$$
\ca{T}^*={1\over(2\pi
i)^2}\int_{(1)}\!\int_{(1)}\prod_p\left(1-{|\Omega^+(p)|-1\over
p-1}\left({1\over p^{s_1}}+{1\over p^{s_2}}-{1\over p^{s_1+s_2}}\right)
\right){R^{s_1+s_2}\over(s_1s_2)^{k+\ell+1}}ds_1ds_2.\leqno(2.11)
$$ 
For $p>H$, we have $|\Omega^+(p)|=k+1$, since $h\notin\ca{H}$.
Thus, we consider the function
$$
\prod_p(\cdots)\left({\zeta(s_1+1)\zeta(s_2+1)\over\zeta(s_1+s_2+1)}
\right)^k\leqno(2.12)
$$
as in $(1.17)$. If $\ca{H}^+$ is admissible, the singular series is
${\eufm S}(\ca{H}^+)$ and the argument and computation of residues is
analogous to above. Thus we find that provided $h\notin\ca{H}$
$$
\ca{T}^*={{\eufm S}(\ca{H}^+)\over(k+2\ell)!}{2\ell\choose\ell}(\log
R)^{k+2\ell}+O((\log N)^{k+2\ell-1}(\log\log N)^c).\leqno(2.13)
$$
On the other hand, if $\ca{H}^+$ is not admissible
or ${\eufm S}(\ca{H}^+)=0$, then the
Euler product in $(2.11)$ vanishes at either $s_1=0$ or $s_2=0$ to the
order equal to the number of primes such that $|\Omega^+(p)|=p$.
However, since we have then $p\le k+1$, the necessary change to the above
reasoning results only in the lack of the main term in $(2.13)$ and
the error term remains to be the same or actually smaller.

Finally, if $h\in\ca{H}$, the above evaluation applies with the translation 
$k\mapsto k-1$, $\ell\mapsto\ell+1$ to (2.13) because of $(2.3)$.
\medskip
From this, we obtain
\medskip
\noindent
{\bf Lemma 2.} {\it Provided $(1.1)$, $(1.2)$, and $(2.4)$ hold, we have,
for $R\le N^{\theta/2}/(\log N)^C$ with a sufficiently
large $C>0$ depending only on $k$ and $\ell$,
$$
\leqalignno{
&\sum_{N<n\le 2N}\varpi(n+h)\Lambda_R(n;\ca{H},k+\ell)^2&(2.14)\cr
=&\cases{\displaystyle{{\eufm S}(\ca{H}\cup\{h\})\over(k+2\ell)!}
{2\ell\choose\ell}N(\log R)^{k+2\ell}+O(N(\log N)^{k+2\ell-1}(\log\log
N)^c)& if $h\not\in\ca{H}$,\cr
\displaystyle{{\eufm S}(\ca{H})\over(k+2\ell+1)!}
{2(\ell+1)\choose\ell+1}N(\log
R)^{k+2\ell+1}+O(N(\log N)^{k+2\ell}(\log\log N)^c)& if $h\in\ca{H}$.
}}
$$
}
\par

\vskip 0.6cm
\centerline{\bf 3. Proof of (0)} 
\medskip
We are now ready to prove $(0)$. We shall evaluate
the expression
$$
\sum_{\scr{\ca{H}\subseteq[1,H]}\atop\scr
{|\ca{H}|=k}}\sum_{N<n\le 2N}\left(\sum_{h\le
H}\varpi(n+h)-\log3N\right)
\Lambda_R(n;\ca{H},k+\ell)^2.\leqno(3.1)
$$
If this turns out to be positive, then there exists an integer
$n\in(N,2N]$ such that
$$
\sum_{h\le H}\varpi(n+h)-\log3N>0.\leqno(3.2)
$$
That is, there exists a subinterval of length $H$ in
$(N,2N+H]$ which contains {\it
two} primes; hence
$$
\min_{N<p_r\le 2N+H}(p_{r+1}-p_r)\le H.\leqno(3.3)
$$
\par
Here, we need to quote, from [2],
$$
\sum_{\scriptstyle \ca{H}\subseteq [1,H]\atop\scriptstyle |\ca{H}|=k}{\eufm 
S}(\ca{H})
=(1+o(1))H^k,\leqno(3.4)
$$
as $H$ tends to infinity. With this and Lemma 1, we see that $(3.1)$ is
asymptotically equal to
$$
\leqalignno{
\sum_{\scr{\ca{H}\subseteq [1,H]}\atop\scr|\ca{H}|=k}&
\sum_{N<n\le 2N}\left\{\sum_{\scr{h\le
H}\atop\scr{h\not\in\ca{H}}}
+\sum_{\scr{h\le H}\atop\scr{h\in\ca{H}}}\right\}\varpi(n+h)
\Lambda_R(n;\ca{H},k+\ell)^2&(3.5)\cr
&-{1\over(k+2\ell)!}{2\ell\choose\ell}NH^k(\log N)(\log R)^{k+2\ell},
}
$$
with an admissible error of the size of $o\left(NH^k(\log
N)^{k+2\ell+1}\right)$. By Lemma 2 and $(3.4)$ with an appropriate
replacement of $\ca{H}$, this is asymptotically equal, in the same sense,
to
$$
\leqalignno{
&{1\over(k+2\ell)!}{2\ell\choose\ell}NH^{k+1}(\log R)^{k+2\ell}
+{k\over(k+2\ell+1)!}{2(\ell+1)\choose\ell+1}
NH^k(\log R)^{k+2\ell+1}&(3.6)\cr
&\qquad-{1\over(k+2\ell)!}{2\ell\choose\ell}NH^k(\log N)(\log
R)^{k+2\ell}\cr
=&\left\{H+{k\over k+2\ell+1}\cdot{2(2\ell+1)\over\ell+1}
\cdot\log R-\log N\right\}
{1\over(k+2\ell)!}{2\ell\choose\ell}NH^k(\log R)^{k+2\ell}. }
$$
Hence $(3.1)$ is positive, provided
$$
H\ge\left(1+\varepsilon-{k\over
k+2\ell+1}\cdot{2(2\ell+1)\over\ell+1}\cdot{\theta\over2}
\right)\log N,\leqno(3.7)
$$
with any fixed $\varepsilon>0$.
Therefore, with $\ell=[\sqrt{k}]$, say, we obtain
$$
\liminf_{n\to\infty}{p_{n+1}-p_n\over\log p_n}\le
\max\left\{0, 1-2\theta\right\}.\leqno(3.8)
$$
In particular, the Bombieri--Vinogradov Prime Number Theorem [1,
Th\'eor\`eme 17] gives rise to the assertion $(0)$. This ends the proof.
\medskip
Finally, we shall exhibit two conjectural assertions: 
\smallskip
\noindent
1) If we have
$\theta>{1\over2}$, then there will be infinitely many
$n$ such that $p_{n+1}-p_n\le c(\theta)$ with an
absolute constant $c(\theta)$. In fact, we would be able to suppose
$H>c(\theta)$ in the above as far as $(3.6)$, and the assertion follows
immediately.
\smallskip
\noindent
2) If we have $\theta>{20\over 21}$, then we will be able to
assert the simultaneous appearance of primes in 
admissible 7-tuples.  For instance, at least 
{\it two} of the seven integers
$\{n,n+2,n+6,n+8,n+12,n+18,n+20\}$ will be primes for infinitely many
$n$. In particular, $p_{n+1}-p_n\le 20$
infinitely often. To prove this, let $\cal H$ be such a tuple, and
consider, in place of $(3.1)$,
$$
\sum_{N<n\le 2N}\left(\sum_{h\in{\cal H}}
\varpi(n+h)-\log3N\right)\Lambda_R(n;\ca{H},8)^2.\leqno(3.9)
$$
Lemmas 1 and 2, with $k=7,\,\ell=1$, imply
that under the present assumption on $\theta$ this is asymptotically equal
to
$$
{2\over9!}{\eufm S}({\cal H})\left({21\over 10}\log R-\log
N\right)N(\log R)^9>0,\leqno(3.10)
$$
provided $R=N^{10/21+\xi}$ with sufficiently large
$N$ and a small $\xi>0$. Hence the assertion follows.

\vskip 0.6cm
\centerline{\bf 4. Expository} 
\medskip
The principal idea in [3] is the amazing effect induced by  the introduction of the
parameter $\ell$ in $(1.9)$. The sieve weight $\mu(d)(\log
m/d)^{k+\ell}$, $d|m$, applied to the polynomial $m=P(n;\ca{H})$ detects
$n$ with which
$P(n;\ca{H})$ has $k+\ell$ distinct prime factors at most, implying
that the integers $n+h_j$, $j\le k$, are mostly primes,
provided $k$ is large compared with $\ell$. By a standard method in this field, we
approximate these sieve weights by
$\lambda_R(d;k+\ell)$, and consider the Selberg sieve
situation $(1.9)$, with the parameters $\ell$ and $R$ at our disposal.
An asymptotic formula for the sum $(1.9)$ is given in $(1.31)$. Then,
to detect at least two primes among $n+h_j$, $j\le k$, a
usual weighted sieve situation is considered at $(3.1)$; for this 
the other asymptotic formula $(2.14)$ is required. The upshot is
condensed in $(3.6)$ and $(3.7)$. The proof of $(0)$ requires that both 
$k$ and $\ell$ can be taken appropriately and the 
Bombieri--Vinogradov prime number theorem is available.
\par
Rendering the above more technically, the reason for success 
lies not only in the introduction of the parameter $\ell$ but
also in the trivial fact $(2.3)$, which   brings forth the translation 
$\ell\mapsto\ell+1$ as remarked in the proof of Lemma 2.
This introduces the factor
$2(\ell+1)\choose\ell+1$ on the right of $(2.14)$. One should note that
${2(\ell+1)\choose\ell+1}/{2\ell\choose\ell}={2(2\ell+1)\over(\ell+1)}$,
which tends to $4$ as $\ell\to\infty$. This is
extremely critical when appealing to the Bombieri--Vinogradov
prime number theorem. On the other hand, the translation $k\mapsto k-1$
does not cause any effect as long as $k$ is much larger than $\ell$.

\vskip 0.6cm
\noindent
\centerline{\bf References}
\medskip
\noindent
\item{[1]} E. Bombieri. Le Grand Crible dans la  
Th\'eorie Analytique des Nombres (Seconde \'Edition). 
Ast\'erisque {\bf 18}, Paris 1987.
\item{[2]} P.X. Gallagher. On the distribution of primes in short
intervals. Mathematika, {\bf 23} (1976), 4--9. Corrigendum, 
ibid {\bf 28} (1981), 86.
\item{[3]} D.A. Goldston, J. Pintz, and C.Y. Y{\i}ld{\i}r{\i}m. Small
gaps between primes II (Preliminary).  Preprint, February 8, 2005.
\item{[4]} E.C. Titchmarsh. The Theory of the Riemann Zeta-Function.
Oxford at the Clarendon Press, First Edition, 1951; Second Edition, 1986.
\vskip 0.8cm
\small
\baselineskip=8pt
\noindent
Department of Mathematics, San Jose State University, San Jose, CA 95192,
USA
\par
\noindent
E-mail address: goldston@math.sjsu.edu
\medskip
\noindent
Department of Mathematics, College of Science and Technology, Nihon
University, Surugadai, Tokyo 101-8308, JAPAN
\par
\noindent
E-mail address: ymoto@math.cst.nihon-u.ac.jp
\medskip
\noindent
R\'enyi Mathematical Institute of the Hungarian Academy of Sciences,
H-1364 Budapest, P.O.B. 127, HUNGARY
\par
\noindent
E-mail address: pintz@renyi.hu
\medskip
\noindent
Department of Mathematics, Bo\~gazi\c{c}i University, Bebek, 
Istanbul 34342 
\par
\noindent
\& Feza G\"ursey Enstit\"us\"u, \c{C}engelk\"oy,
Istanbul, P.K. 6, 81220 TURKEY
\par
\noindent
E-mail address: yalciny@boun.edu.tr
\medskip
\noindent

\bye